\documentclass[10pt]{article}

\usepackage{amsmath}
\usepackage{amssymb}
\usepackage{indentfirst}
\usepackage{graphics} 
\usepackage{color}

\makeatletter

\@addtoreset{equation}{section}
\makeatother
\def\<{\langle}
\def\>{\rangle}

\newtheorem{lem}{Lemma}[section]
\newtheorem{theo}{Theorem}[section]
\newtheorem{rem}{Remark}[section]
\newtheorem{pro}{Proposition}[section]

\makeatletter
   
   \@addtoreset{equation}{section}
\makeatother

\begin{document}
\title{\bf $L^{2}$-blowup estimates of the wave equation and its application to local energy decay}

\author{Ryo Ikehata\thanks{ikehatar@hiroshima-u.ac.jp} \\ {\small Department of Mathematics}\\ {\small Division of Educational Sciences}\\ {\small Graduate School of Humanities and Social Sciences} \\ {\small Hiroshima University} \\ {\small Higashi-Hiroshima 739-8524, Japan}}
\date{}
\maketitle
\begin{abstract}
We consider the Cauchy problems in ${\bf R}^{n}$ for the wave equation with a weighted $L^{1}$-initial data. We derive sharp infinite time blowup estimates of the $L^{2}$-norm of solutions in the case of $n = 1$ and $n = 2$. Then, we apply it to the local energy decay estimates for $n = 2$, which is not studied so completely when the $0$-th moment of the initial velocity does not vanish. 
The idea to derive them is strongly inspired from a technique used in \cite{I-14, IO}.  
\end{abstract}
\section{Introduction}
\footnote[0]{Keywords and Phrases: Wave equation; weighted $L^{1}$-data; local energy; low dimensional case; blowup in infinite time; sharp estimates.}
\footnote[0]{2010 Mathematics Subject Classification. Primary 35L05; Secondary 35B40, 35C20.}

We consider the Cauchy problem of the wave equation:
\begin{align}
& u_{tt} -\Delta u = 0,\ \ \ (t,x)\in (0,\infty)\times {\bf R}^{n},\label{eqn}\\
& u(0,x)= u_0(x), \quad  u_{t}(0,x)= u_{1}(x),\ x\in{\bf R}^{n}.\label{initial}
\end{align}
Here, we assume, for the moment, $[u_{0},u_{1}] \in H^{1}({\bf R}^{n})\times L^{2}({\bf R}^{n})$.\\ 
\noindent
Concerning the existence of a unique energy solution to problem \eqref{eqn}-\eqref{initial}, by the Lumer-Phillips Theorem one can find that the problem (1.1)-(1.2) has a unique weak solution
\[u \in C([0,\infty);H^{1}({\bf R}^{n})) \cap C^{1}([0,\infty);L^{2}({\bf R}^{n})) \]
satisfying the energy conservation law such that
\begin{equation}\label{i-1}
E(t) = E(0),\quad t \geq 0,
\end{equation}
where the total energy $E(t)$ for the solution to problem \eqref{eqn}-\eqref{initial} can be defined by
\[E(t) := \frac{1}{2}\left(\Vert u_{t}(t,\cdot)\Vert^{2} + \Vert\nabla u(t,\cdot)\Vert^{2}\right).\] 
Here, $\Vert u\Vert$ denotes the usual $L^{2}$-norm of $u \in L^{2}({\bf R}^{n})$.

Our main purpose is to establish the sharp estimates in time of the factor $M(t) := \Vert u(t,\cdot)\Vert$. This observation is important from the view point of the local energy decay, which is one of main topics in the wave equation field. As is known from the work of C. Morawetz \cite{M}, the main task to get the local energy decay is to control the behavior of the function $M(t)$.  By the way, as in a series of papers \cite{AI, IN} and so on, in the exterior mixed problems on a star-shaped compliment (or more generally just on exterior domains) one can control $M(t)$ for all $n \geq 2$. While if one considers the local energy decay problems to (1.1) in the Euclidean space ${\bf R}^{n}$ (no obstacle case!) one can proceed the same argument as in the exterior mixed problem case when $n \geq 3$ (see \cite{CI}), however, in the case of $n = 2$, one has to impose stronger assumption such that $I_{0} := \int_{{\bf R}^{2}}u_{1}(x)dx = 0$, which means the vanishing condition of the $0$th moment of the initial data. Note that in these problems, the author images that one never assumes the compactness of the support of the initial data, and low regularity on them to control the function $M(t)$, and to get the local energy decay. Therefore, it is highly desirable to remove the stronger condition $I_{0} = 0$ in the two dimensional case for the completeness of the theory. In this connection, for $n \geq 3$ $L^{p}$-decay problems of the solution to problem (1.1)-(1.2) including the Klein-Gordon type are already discussed by \cite{Strauss} and \cite{W}, however, there the compactness of the support of the initial data seems to be conscious, and their main concern is not in local energy decay itself, but an application to nonlinear problems of the $L^{p}$-boundedness.\\

Now, in this paper we study the asymptotic behavior of the function $M(t)$ as $t \to \infty$ under the non-compactly supported initial data framework. It seems that this problem does not studied so far. We employ the simple Fourier analysis to get the estimates of the function $M(t)$ by dividing the Fourier space of the integrand of the Fourier transformed solution $\hat{u}(t,\xi)$ into two parts: one is low frequency zone ($\vert\xi\vert \leq q(t)$), and the other is high frequency zone ($\vert\xi\vert \geq q(t)$), where $q(t)$ is a $t$-variable function chosen suitably. So, one never relies on the Hardy type inequality, which never holds in the low dimensional case (cf. \cite{O}) because a core of the idea in a series of papers \cite{AI, CI, I-1,I-2, IN, IS} is that the boundedness of $M(t)$ can be derived by using the Hardy inequality. The problem to study the behavior of $M(t)$ as $t \to \infty$ is important from the view point that the Hardy type inequality does not hold in the low dimensional case. \\    
 
Before going to introduce our Theorems, we set
\[I_{0,n} := \Vert u_{0}\Vert + \Vert u_{0}\Vert_{L^{1}({\bf R}^{n})} +\Vert u_{1}\Vert + \Vert u_{1}\Vert_{L^{1}({\bf R}^{n})}.\]
Our first result is concerned with the optimal infinite time blowup result in the case of $n = 1$.
\begin{theo}\label{theorem1}
Let $n = 1$. Let $[u_{0},u_{1}] \in H^{1}({\bf R}) \times L^{2}({\bf R})$. Then, the solution $u(t,x)$ to problem \eqref{eqn}-\eqref{initial} satisfies the following properties under the additional regularity on the initial data:
\[[u_{0},u_{1}] \in L^{1}({\bf R}) \times L^{1}({\bf R}) \quad \Rightarrow \quad \Vert u(t,\cdot)\Vert \leq C_{1}I_{0,1}\sqrt{t},\]
\[[u_{0},u_{1}] \in L^{1}({\bf R}) \times L^{1,1}({\bf R}) \quad \Rightarrow \quad C_{2}\left\vert\int_{{\bf R}}u_{1}(x)dx\right\vert\sqrt{t} \leq \Vert u(t,\cdot)\Vert\]
for $t \gg 1$, where $C_{j} > 0$ {\rm ($j = 1,2$)} are constants depending only on the space dimension.
\end{theo}
Our next result is the case of $n = 2$.
\begin{theo}\label{theorem2}
Let $n = 2$. Let $[u_{0},u_{1}] \in H^{1}({\bf R}^{2}) \times L^{2}({\bf R}^{2})$. Then, the solution $u(t,x)$ to problem \eqref{eqn}-\eqref{initial} satisfies the following properties under the additional regularity on the initial data:
\[[u_{0},u_{1}] \in L^{1}({\bf R}^{2}) \times L^{1}({\bf R}^{2}) \quad \Rightarrow \quad \Vert u(t,\cdot)\Vert \leq C_{1}I_{0,2}\sqrt{\log t},\]
\[[u_{0},u_{1}] \in L^{1}({\bf R}^{2}) \times L^{1,1}({\bf R}^{2}) \quad \Rightarrow \quad C_{2}\left\vert\int_{{\bf R}^{2}}u_{1}(x)dx\right\vert\sqrt{\log t} \leq \Vert u(t,\cdot)\Vert\]
for $t \gg 1$, where $C_{j} > 0$ {\rm ($j = 1,2$)} are constants depending only on the space dimension. 
\end{theo}
\begin{rem}{\rm Surprisingly if one imposes the additional $L^{1,1}$-regularity on the initial velocity, one can get the sharp infinite time blowup result in the one and two dimensional free waves. It seems that this crucial phenomenon is not so known. Generally, people assume a vanishing moment condition $\int_{{\bf R}^{n}}u_{1}(x)dx\ = 0$ when one studies the $L^{2}$-boundedness of the solution $u(t,x)$, for the low dimensional case $n = 1,2$ (cf., \cite{Azer-ike}). 
}
\end{rem}
\begin{rem}{\rm From the proof done in Section 3 for the upper bound estimates one can see that $L^{1,1}$-assumption is no need to get upper bound estimates, and only $L^{1}$-assumption works well. On the other hand, to obtain lower bound estimates one needs to impose the $L^{1,1}$ assumption on the initial velocity $u_{1}$ (see \eqref{i-13} and \eqref{i-13-1}).}
\end{rem}
\noindent
Let us introduce one example to make sure the reliability of Theorems in the paper about only in the one dimensional case.\\
{\bf Example.} Let $n = 1$, $u_{0} \equiv 0$, and we choose $u_{1} \in (L^{1,1}({\bf R})\cap L^{2}({\bf R}))$ as follows:\\
\[ u_{1}(x) = \left\{
     \begin{array}{ll}
       \displaystyle{2}&
            \qquad (\vert x\vert \leq 1), \\[0.2cm]
        \displaystyle{0}& \qquad (\vert x\vert > 1).
     \end{array} \right. \]  
Then, one easily gets the formula:
\[u(t,x) = \frac{1}{2}\int_{x-t}^{x+t}u_{1}(r)dr = -\chi(x-t) + \chi(x+t),\]
where
\[\chi(x) := \frac{1}{2}\int_{0}^{x}u_{1}(r)dr.\]
One can soon check that
\[ -\chi(x-t) = \left\{
     \begin{array}{ll}
       \displaystyle{-(x-t)}&
            \qquad (\vert x-t\vert \leq 1), \\[0.2cm]
      \displaystyle{-1}&
            \qquad (x-t > 1), \\[0.2cm]       
       \displaystyle{1}& \qquad (x-t < -1),
     \end{array} \right. \]  
\[ \chi(x+t) = \left\{
     \begin{array}{ll}
       \displaystyle{(x+t)}&
            \qquad (\vert x+t\vert \leq 1), \\[0.2cm]
      \displaystyle{1}&
            \qquad (x+t > 1), \\[0.2cm]       
       \displaystyle{-1}& \qquad (x+t < -1).
     \end{array} \right. \]  
Therefore, if one fixes a sufficiently large $t_{0} > 1$, then one can finally arrive at the following situation for the solution $u(t_{0},x)$:  
\[ u(t_{0},x) = \left\{
     \begin{array}{ll}
       \displaystyle{2}&
            \qquad (1-t_{0} \leq x \leq t_{0}-1), \\[0.2cm]
      \displaystyle{x + t_{0} + 1}&
            \qquad (-t_{0}-1 < x < 1-t_{0}), \\[0.2cm]       
       \displaystyle{-x + t_{0} + 1}&
            \qquad (t_{0}-1 < x < t_{0}+1). \\[0.2cm]
        \displaystyle{0}& \qquad (x < -t_{0}-1,\, t_{0}+1 < x).    
     \end{array} \right. \]  
Therefore, one can compute
\[\int_{{\bf R}}\vert u(t_{0},x)\vert^{2}dx = 8(t_{0}-1) + \frac{16}{3},\]
which implies the credibility of the result
\[\Vert u(t_{0},\cdot)\Vert \sim \sqrt{t_{0}}, \quad t_{0} \gg 1.\]
In this connection, one naturally has the condition $\displaystyle{\int_{{\bf R}}}u_{1}(x)dx = 4 \ne 0$. This example can be stated in the fundamental text book of PDE (e.g. \cite{K}). 
\begin{rem}{\rm In a sense, the infinite time blowup results in $L^{2}$-sense may express quantitatively failure of the Huygens principle.}
\end{rem}
{\bf Notation.} {\small Throughout this paper, $\| \cdot\|_q$ stands for the usual $L^q({\bf R}^{n})$-norm. For simplicity of notation, in particular, we use $\| \cdot\|$ instead of $\| \cdot\|_2$. 
We also introduce the following weighted functional spaces.
\[L^{1,\gamma}({\bf R}^{n}) := \left\{f \in L^{1}({\bf R}^{n}) \; \bigm| \; \Vert f\Vert_{1,\gamma} := \int_{{\bf R}^{n}}(1+\vert x\vert^{\gamma})\vert f(x)\vert dx < +\infty\right\}.\]
One denotes the Fourier transform ${\cal F}_{x\to\xi}(f)(\xi)$ of $f(x)$ by 
\[{\cal F}_{x\to\xi}(f)(\xi) = \hat{f}(\xi) := \displaystyle{\int_{{\bf R}^{n}}}e^{-ix\cdot\xi}f(x)dx, \quad \xi \in {\bf R}^n,\]
as usual with $i := \sqrt{-1}$, and ${\cal F}_{\xi\to x}^{-1}$ expresses its inverse Fourier transform. Finally, we denote the surface area of the $n$-dimensional unit ball by $\omega_{n} := \displaystyle{\int_{\vert\omega\vert = 1}}d\omega$. }

The paper is organized as follows. In Section 2 we derive the lower bound estimates of the $L^{2}$-norm of solutions, and in Section 3 we obtain the upper bound estimates of the $L^{2}$-norm of solutions, and by combining the results obtained in Sections $2$ and $3$ one can prove Theorems 1.1 and 1.2 at a stroke. Section 4 is devoted to apply Theorem 1.2 to the local energy decay results of the wave equation (1.1) in the two dimensional Euclidean space ${\bf R}^{2}$. \\


\section{$L^{2}$-lower bound estimates of the solution}

In this section, let us derive the lower bound estimates of the $t$-function $\Vert u(t,\cdot)\Vert $ by using the Plancherel Theorem and low frequency estimates. This technique is well-developed in the damped wave equation field (cf. \cite{I-14}).\\ 
We first prepare two basic facts which will be frequently used in the proof.\\
Set
\begin{equation}\label{i-2}
L := \sup_{\theta \ne 0}\left\vert \frac{\sin\theta}{\theta}\right\vert < +\infty.
\end{equation}
Furthermore, since
\[\lim_{\theta \to +0}\frac{\sin\theta}{\theta} = 1,\]
one can find a real number $\delta_{0} \in (0,1)$ such that
\begin{equation}\label{i-3}
\left\vert \frac{\sin\theta}{\theta}\right\vert \geq \frac{1}{2}
\end{equation}
for all $\theta \in (0,\delta_{0}]$. One also prepares the fundamental inequality:
\begin{equation}\label{i-4}
\vert a + b\vert^{2} \geq \frac{1}{2}\vert a\vert^{2} - \vert b\vert^{2} 
\end{equation}
for all $a, b \in {\bf C}$.\\
In order to get the lower bound estimate for the quantity $\Vert u(t,\cdot)\Vert$, it suffices to treat $\Vert w(t,\cdot)\Vert$ with $w(t,\xi) := {\cal F}_{x \to \xi}(u(t,\cdot))(\xi) = \hat{u}(t,\xi)$ because of the Plancherel Theorem.\\
Now we decompose the quantity $\Vert w(t,\cdot)\Vert$ as follows: for each $n \geq 1$
\begin{equation}\label{i-5}
\Vert w(t,\cdot)\Vert^{2} = \left(\int_{\vert\xi\vert \leq \frac{\delta_{0}}{t}} + \int_{\vert\xi\vert \geq \frac{\delta_{0}}{t}}\right)\vert w(t,\xi)\vert^{2}d\xi = I_{low}^{(n)}(t) + I_{high}^{(n)}(t).
\end{equation}
Here we have just chosen $t > 0$ large enough such that
\[\frac{\delta_{0}}{t} \leq 1.\]

By the way, as is well-known in the text book of PDE, in the Fourier space ${\bf R}_{\xi}^{n}$ the problem \eqref{eqn}-\eqref{initial} and its solution $u(t,x)$ can be transformed into the following ODE with parameter $\xi \in {\bf R}_{\xi}^{n}$:
\begin{align}
& w_{tt} + |\xi|^2 w = 0,\ \ \ t>0,\quad \xi \in {\bf R}_{\xi}^{n},\label{eqnfourier}\\
& w(0,\xi)= w_0(\xi), \quad  w_{t}(0,\xi)= w_{1}(\xi),\ \ \ \xi \in{\bf R}^{n} ,\label{initialfourier}
\end{align}
where $w_{0}(\xi) := \hat{u}_0(\xi)$ and $w_{1}(\xi) := \hat{u}_1(\xi)$. Moreover, one can easily solve the problem \eqref{eqnfourier}-\eqref{initialfourier} as follows:
\begin{equation}\label{i-6}
w(t,\xi) = \frac{\sin(t\vert\xi\vert)}{\vert\xi\vert}w_{1}(\xi) + \cos(t\vert\xi\vert)w_{0}(\xi).
\end{equation}

Now, let us give the lower bound estimates for the quantity $I_{low}^{(n)}(t)$ because of $\Vert w(t,\cdot)\Vert^{2} \geq I_{low}^{(n)}(t)$. The technique developed below is strongly inspired from the idea in \cite{I-14}.  Indeed, it follows from \eqref{i-6} and \eqref{i-4} that
\[I_{low}^{(n)}(t) = \int_{\vert\xi\vert \leq \frac{\delta_{0}}{t}}\left\vert\frac{\sin(t\vert\xi\vert)}{\vert\xi\vert}w_{1}(\xi) + \cos(t\vert\xi\vert)w_{0}(\xi)\right\vert^{2} d\xi\]
\[\geq \frac{1}{2}\int_{\vert\xi\vert \leq \frac{\delta_{0}}{t}}\frac{\sin^{2}(t\vert\xi\vert)}{\vert\xi\vert^{2}}\vert w_{1}(\xi)\vert^{2}d\xi -  \int_{\vert\xi\vert \leq \frac{\delta_{0}}{t}}\cos^{2}(t\vert\xi\vert)\vert w_{0}(\xi)\vert^{2}d\xi\]
\begin{equation}\label{i-7}
=: \frac{1}{2}J_{1}(t) - J_{2}(t).
\end{equation}
Let us first estimate $J_{1}(t)$ by using the decomposition of the initial data $w_{1}(\xi)$ in the Fourier space: 
\[w_{1}(\xi) = P + (A(\xi)-iB(\xi)),\quad \xi \in {\bf R}_{\xi}^{n},\]
where
\[P := \int_{{\bf R}^{n}}u_{1}(x)dx,\]
\[A(\xi) := \int_{{\bf R}^{n}}(\cos(x\xi)-1)u_{1}(x)dx, \quad B(\xi) := \int_{{\bf R}^{n}}\sin(x\xi)u_{1}(x)dx.\]
It is known (see \cite{I-04}) that with some constant $M > 0$ one has
\begin{equation}\label{i-8}
\vert A(\xi)-iB(\xi)\vert \leq M\vert\xi\vert\Vert u_{1}\Vert_{L^{1,1}},\quad \xi \in {\bf R}_{\xi}^{n},
\end{equation}
in the case when $u_{1} \in L^{1,1}({\bf R}^{n})$. Then, it follows from \eqref{i-4} that
\[J_{1}(t) = \int_{\vert\xi\vert \leq \frac{\delta_{0}}{t}}\frac{\sin^{2}(t\vert\xi\vert)}{\vert\xi\vert^{2}}\vert w_{1}(\xi)\vert^{2}d\xi\]
\[\geq \frac{P^{2}}{2}\int_{\vert\xi\vert \leq \frac{\delta_{0}}{t}}\frac{\sin^{2}(t\vert\xi\vert)}{\vert\xi\vert^{2}}d\xi - \int_{\vert\xi\vert \leq \frac{\delta_{0}}{t}}\vert A(\xi)-iB(\xi)\vert^{2}\frac{\sin^{2}(t\vert\xi\vert)}{\vert\xi\vert^{2}}d\xi\]
\begin{equation}\label{i-9}
= \frac{P^{2}}{2}K_{1}(t)-K_{2}(t).
\end{equation}
$K_{2}(t)$ can be estimated from above as follows by using \eqref{i-8}:
\[K_{2}(t) \leq M^{2}\Vert u_{1}\Vert_{L^{1,1}}^{2}\int_{\vert\xi\vert \leq \frac{\delta_{0}}{t}}\sin^{2}(t\vert\xi\vert)d\xi = M^{2}\Vert u_{1}\Vert_{L^{1,1}}^{2}\omega_{n}\int_{0}^{\frac{\delta_{0}}{t}}r^{n-1}dr\]
\begin{equation}\label{i-10}
= \frac{M^{2}}{n}\omega_{n}\delta_{0}^{n}\Vert u_{1}\Vert_{L^{1,1}}^{2} t^{-n},\quad t \gg 1.
\end{equation}
While, one can get the lower bound estimate for $K_{1}(t)$ because of \eqref{i-3}:
\[K_{1}(t) = t^{2}\int_{\vert\xi\vert \leq \frac{\delta_{0}}{t}}\frac{\sin^{2}(t\vert\xi\vert)}{(t\vert\xi\vert)^{2}}d\xi \geq \frac{t^{2}}{4}\int_{\vert\xi\vert \leq \frac{\delta_{0}}{t}}d\xi \]
\begin{equation}\label{i-11}
= \frac{\omega_{n}\delta_{0}^{n}}{4n}t^{2-n}, \quad t \gg 1.
\end{equation}
Therefore from \eqref{i-9}, \eqref{i-10} and \eqref{i-11} one can get the estimate from below for $J_{1}(t)$:
\begin{equation}\label{i-12}
J_{1}(t) \geq \frac{P^{2}}{2}\frac{\omega_{n}\delta_{0}^{n}}{4n}t^{2-n}-\frac{M^{2}}{n}\omega_{n}\delta_{0}^{n}\Vert u_{1}\Vert_{L^{1,1}}^{2} t^{-n} \quad t \gg 1.
\end{equation}
Since the upper bound estimate of $J_{2}(t)$ can be easily obtained as follows: 
\[J_{2}(t) \leq \int_{\vert\xi\vert \leq \frac{\delta_{0}}{t}}\vert w_{0}(\xi)\vert^{2}d\xi \leq \Vert u_{0}\Vert_{L^{1}}^{2}\omega_{n}\int_{0}^{\frac{\delta_{0}}{t}}r^{n-1}dr = \frac{\omega_{n}\delta_{0}^{n}}{n}\Vert u_{0}\Vert_{L^{1}}^{2}t^{-n},\]
by \eqref{i-12}, \eqref{i-5} and \eqref{i-7} one has just arrived at the following crucial lower bound estimate for $\Vert w(t,\cdot)\Vert$:
\begin{equation}\label{i-13}
\Vert w(t,\cdot)\Vert^{2} \geq I_{low}^{n}(t) \geq  \frac{P^{2}}{4}\frac{\omega_{n}\delta_{0}^{n}}{4n}t^{2-n}-\frac{M^{2}}{n}\omega_{n}\delta_{0}^{n}\Vert u_{1}\Vert_{L^{1,1}}^{2} t^{-n} -\frac{\omega_{n}\delta_{0}^{n}}{n}\Vert u_{0}\Vert_{L^{1}}^{2}t^{-n}\quad t \gg 1.
\end{equation}
Thus, there exists a positive real number $t_{0}$ such that 
\begin{equation}\label{i-13-1}
\Vert w(t,\cdot)\Vert^{2} \geq \frac{P^{2}}{32n}\omega_{n}\delta_{0}^{n} t^{2-n}
\end{equation}
for all $t \geq t_{0}$. It should be mentioned that $t_{0} > 0$ depends on $n$ and the quantities $\Vert u_{1}\Vert_{L^{1.1}}$ and $\Vert u_{0}\Vert_{L^{1}}$.\\
By choosing $n = 1$ in \eqref{i-13-1} one has the following blowup property of the solution when $P \ne 0$.
\begin{lem}\label{lem1}Let $n = 1$, and $[u_{0},u_{1}] \in L^{1}({\bf R}^{n}) \times L^{1,1}({\bf R}^{n})$. Then, it holds that
\[\Vert w(t,\cdot)\Vert^{2} \geq CP^{2}t, \quad t \gg 1.\]
\end{lem}

Let us prove the statement for $n = 2$ at a stroke by using a trick with a function $e^{-r^{2}}$. Indeed, it follows from \eqref{i-4} and a similar argument to the one dimensional case that 
\[\Vert w(t,\cdot)\Vert^{2} \geq \frac{1}{2}\int_{{\bf R}^{2}}\frac{\sin^{2}(tr)}{r^{2}}\vert w_{1}(\xi)\vert^{2}d\xi -\int_{{\bf R}^{2}}\cos^{2}(tr)\vert w_{0}(\xi)\vert^{2}d\xi\]
\[\geq \frac{1}{2}\int_{{\bf R}^{2}}e^{-r^{2}}\frac{\sin^{2}(tr)}{r^{2}}e^{r^{2}}\vert P+(A(\xi)-iB(\xi))\vert^{2}d\xi- \Vert u_{0}\Vert^{2}\]
\[\geq \frac{1}{4}P^{2}\int_{{\bf R}^{2}}e^{-r^{2}}\frac{\sin^{2}(tr)}{r^{2}}d\xi - \frac{1}{2}\int_{{\bf R}^{2}}e^{-r^{2}}\frac{\sin^{2}(tr)}{r^{2}}\left(M^{2}\Vert u_{1}\Vert_{L^{1,1}}^{2}r^{2}\right)d\xi- \Vert u_{0}\Vert^{2}\]
\[\geq \frac{1}{4}P^{2}\int_{{\bf R}^{2}}e^{-r^{2}}\frac{\sin^{2}(tr)}{r^{2}}d\xi - \frac{1}{2}M^{2}\Vert u_{1}\Vert_{L^{1,1}}^{2}\int_{{\bf R}^{2}}e^{-r^{2}}d\xi -  \Vert u_{0}\Vert^{2}\]
\begin{equation}\label{ike-51}
=:  \frac{1}{4}P^{2}T(t) -\frac{\omega_{2}}{4}M^{2}\Vert u_{1}\Vert_{L^{1,1}}^{2} - \Vert u_{0}\Vert^{2}.
\end{equation}
Now, we apply an idea which has its origin in \cite{IO}. For this purpose, we set
\[\nu_{j} := (\frac{1}{4}+j)\frac{\pi}{t}, \quad \mu_{j} := (\frac{3}{4}+j)\frac{\pi}{t}\quad (j = 0,1,2,\cdots),\]
and choose $t > 1$ large enough such that $\nu_{1} = \frac{5\pi}{4t} < 1$. Then, since $$\vert\sin(tr)\vert \geq \frac{1}{\sqrt{2}}$$
for $r \in [\nu_{j},\mu_{j}]$ and for each $j = 0,1,2,\cdots$, one has the estimate:
\begin{equation}\label{ike-53}
T(t) = \int_{{\bf R}^{2}}e^{-r^{2}}\frac{\sin^{2}(tr)}{r^{2}}d\xi \geq \frac{1}{2}\sum_{j= 0}^{\infty}\int_{\nu_{j} \leq r \leq \mu_{j}}\frac{e^{-r^{2}}}{r^{2}}d\xi = \frac{\omega_{2}}{2}\left(\sum_{j= 0}^{\infty}\int_{\nu_{j}}^{\mu_{j}}\frac{e^{-r^{2}}}{r}dr\right)
\end{equation}
\begin{equation}\label{ike-54}
\geq \frac{\omega_{2}}{2}\left(\frac{1}{2}\int_{\nu_{1}}^{\infty}\frac{e^{-r^{2}}}{r}dr\right) \geq \frac{\omega_{2}}{4}\int_{\nu_{1}}^{1}\frac{e^{-r^{2}}}{r}dr
\end{equation}
\begin{equation}\label{ike-52}
\geq \omega_{2}\frac{e^{-1}}{4}\int_{\nu_{1}}^{1}\frac{1}{r}dr = \omega_{2}\frac{e^{-1}}{4}(\log t +\log 4 -\log 5\pi),
\end{equation}
where in the inequality from \eqref{ike-53} to \eqref{ike-54} one has just used the monotone decreasing property of the function $r \mapsto \frac{e^{-r^{2}}}{r}$, and the fact that $\frac{\pi}{2t} = \mu_{j}-\nu_{j} = \nu_{j+1} -\mu_{j}$ ($\forall j$). Therefore, by \eqref{ike-51} and \eqref{ike-52} one has the following estimates for $n = 2$ for large $t > 1$.
\begin{lem}\label{lem2}Let $n = 2$, and $[u_{0},u_{1}] \in L^{1}({\bf R}^{n}) \times L^{1,1}({\bf R}^{n})$. Then, it holds that
\[\Vert w(t,\cdot)\Vert^{2} \geq CP^{2}\log t, \quad t \gg 1.\]
\end{lem}


\section{$L^{2}$-upper bound estimates of the solution}

In this section, one derives upper bound estimates of the quantity $\Vert u(t,\cdot)\Vert$ as $t \to \infty$ by calculating the function $w(t,\xi)$ in both high and low frequency region.\\ 

From \eqref{i-5} one first shows the upper bound estimate for $I_{low}^{(n)}(t)$ for all $n \geq 1$. Indeed, by \eqref{i-7} one has
\[I_{low}^{(n)}(t) \leq \int_{\vert\xi\vert \leq \frac{\delta_{0}}{t}}\vert\frac{\sin(t\vert\xi\vert)}{\vert\xi\vert}w_{1}(\xi) + \cos(t\vert\xi\vert)w_{0}(\xi) \vert^{2}d\xi\]
\[\leq 2\int_{\vert\xi\vert \leq \frac{\delta_{0}}{t}}\frac{\sin^{2}(t\vert\xi\vert)}{\vert\xi\vert^{2}}\vert w_{1}(\xi)\vert^{2}d\xi  + 2\int_{\vert\xi\vert \leq \frac{\delta_{0}}{t}}\cos^{2}(t\vert\xi\vert)\vert w_{0}(\xi)\vert^{2}d\xi\]
\begin{equation}\label{i-14}
= 2L_{1}(t) + 2L_{2}(t). 
\end{equation}
To begin with, it is easy to get the estimate for $L_{2}(t)$ as in the estimate for $J_{2}(t)$ in Section 2: 
\[L_{2}(t) \leq \int_{\vert\xi\vert \leq \frac{\delta_{0}}{t}}\vert w_{0}(\xi)\vert^{2}d\xi \leq \Vert u_{0}\Vert_{L^{1}}^{2}\omega_{n}\int_{0}^{\frac{\delta_{0}}{t}}r^{n-1}dr\]
\begin{equation}\label{i-15}
 = \frac{\omega_{n}\delta_{0}^{n}}{n}\Vert u_{0}\Vert_{L^{1}}^{2}t^{-n}.
\end{equation}
Let us estimate $L_{1}(t)$ as follows: 
\[L_{1}(t) = \int_{\vert\xi\vert \leq \frac{\delta_{0}}{t}}\frac{\sin^{2}(t\vert\xi\vert)}{\vert\xi\vert^{2}}\vert w_{1}(\xi)\vert^{2}d\xi \]
\begin{equation}\label{i-15-1}
\leq \Vert u_{1}\Vert_{L^{1}}^{2}\int_{\vert\xi\vert \leq \frac{\delta_{0}}{t}}\frac{\sin^{2}(t\vert\xi\vert)}{\vert\xi\vert^{2}}d\xi.
\end{equation}
As in \eqref{i-15} one has
\begin{equation}\label{i-17}
\int_{\vert\xi\vert \leq  \frac{\delta_{0}}{t}}d\xi \leq \frac{\omega_{n}\delta_{0}^{n}}{n}t^{-n},\quad t \gg 1.
\end{equation}
Therefore, if one uses \eqref{i-2}, \eqref{i-15-1} and \eqref{i-17} then it holds that
\begin{equation}\label{i-16}
L_{1}(t) \leq \Vert u_{1}\Vert_{L^{1}}^{2}L^{2}t^{2} \frac{\omega_{n}\delta_{0}^{n}}{n}t^{-n} = \Vert u_{1}\Vert_{L^{1}}^{2}L^{2}\frac{\omega_{n}\delta_{0}^{n}}{n}t^{2-n},\quad t \gg 1.
\end{equation}
\noindent
Thus, from \eqref{i-14}, \eqref{i-15} and \eqref{i-16} one can obtain the low-frequency estimate
\begin{equation}\label{i-low}
I_{low}^{(n)}(t) \leq C\left(\Vert u_{0}\Vert_{L^{1}}^{2}t^{-n} + \Vert u_{1}\Vert_{L^{1}}^{2}t^{2-n}\right), \quad t \gg 1,
\end{equation}
with some constant $C > 0$. 

Next, one treats $I_{high}^{(n)}(t)$ to get the upper bound estimate. Indeed, similar to the above estimate one stars with the following inequalities.
\[I_{high}^{(n)}(t) \leq \int_{\vert\xi\vert \geq \frac{\delta_{0}}{t}}\vert\frac{\sin(t\vert\xi\vert)}{\vert\xi\vert}w_{1}(\xi) + \cos(t\vert\xi\vert)w_{0}(\xi) \vert^{2}d\xi\]
\[\leq 2\int_{\vert\xi\vert \geq \frac{\delta_{0}}{t}}\frac{\sin^{2}(t\vert\xi\vert)}{\vert\xi\vert^{2}}\vert w_{1}(\xi)\vert^{2}d\xi  + 2\int_{\vert\xi\vert \geq \frac{\delta_{0}}{t}}\cos^{2}(t\vert\xi\vert)\vert w_{0}(\xi)\vert^{2}d\xi\]
\begin{equation}\label{i-20}
= 2N_{1}^{(n)}(t) + 2N_{2}^{(n)}(t). 
\end{equation}
It is easy to treat $N_{2}^{(n)}(t)$ as follows. This can hold for all $n \geq 1$. 
\begin{equation}\label{i-21}
N_{2}^{(n)}(t) \leq \int_{\vert\xi\vert \geq \frac{\delta_{0}}{t}}\vert w_{0}(\xi)\vert^{2}d\xi \leq \Vert u_{0}\Vert^{2}.
\end{equation}

Let us estimate $N_{1}^{(n)}(t)$ when $n = 1,2$. First, for all $n \geq 1$ one has 
\[N_{1}^{(n)}(t) = \int_{\vert\xi\vert \geq \frac{\delta_{0}}{\sqrt{t}}}\frac{\sin^{2}(t\vert\xi\vert)}{\vert\xi\vert^{2}}\vert w_{1}(\xi)\vert^{2}d\xi + \int_{\frac{\delta_{0}}{t} \leq \vert\xi\vert \leq \frac{\delta_{0}}{\sqrt{t}}}\frac{\sin^{2}(t\vert\xi\vert)}{\vert\xi\vert^{2}}\vert w_{1}(\xi)\vert^{2}d\xi \]
\begin{equation}\label{i-22}
=: O_{1}^{(n)}(t) + O_{2}^{(n)}(t).
\end{equation}
To begin with, we consider the case $n=1$. In this case, for $O_{1}^{(1)}(t)$ one can get the estimate
\[O_{1}^{(1)}(t) \leq \frac{t}{\delta_{0}^{2}}\int_{\vert\xi\vert \geq \frac{\delta_{0}}{\sqrt{t}}}\sin^{2}(t\vert\xi\vert)\vert w_{1}(\xi)\vert^{2}d\xi\]
\begin{equation}\label{i-23}
\leq \frac{t}{\delta_{0}^{2}}\int_{\vert\xi\vert \geq \frac{\delta_{0}}{\sqrt{t}}}\vert w_{1}(\xi)\vert^{2}d\xi \leq  \frac{t}{\delta_{0}^{2}}\Vert u_{1}\Vert^{2}, \quad t \gg 1.
\end{equation}
\noindent
While, $O_{2}^{(1)}(t)$ can be treated as follows:
\[O_{2}^{(1)}(t) \leq \Vert u_{1}\Vert_{1}^{2}\int_{\frac{\delta_{0}}{t} \leq \vert\xi\vert \leq \frac{\delta_{0}}{\sqrt{t}}}\frac{1}{\vert\xi\vert^{2}}d\xi\]
\begin{equation}\label{i-25}
= \Vert u_{1}\Vert_{1}^{2}\omega_{1}\int_{\frac{\delta_{0}}{t}}^{\frac{\delta_{0}}{\sqrt{t}}}\frac{1}{r^{2}}dr = \frac{\omega_{1}}{\delta_{0}}\Vert u_{1}\Vert_{1}^{2}(t-\sqrt{t}), \quad t \gg 1.
\end{equation}
Thus, from \eqref{i-20}, \eqref{i-21}, \eqref{i-22}, \eqref{i-23} and \eqref{i-25} one has the estimate for $I_{high}^{(n)}(t)$ in the case of $n = 1$:
\begin{equation}\label{i-24}
I_{high}^{(1)}(t)= C\left(\Vert u_{0}\Vert^{2} + \Vert u_{1}\Vert^{2}t + \Vert u_{1}\Vert_{L^{1}}^{2}(t-\sqrt{t})\right), \quad t \gg 1.
\end{equation}

Next, let us give sharp estimates for $N_{1}^{(2)}(t)$ in the case $n = 2$ at a stroke by using a more precise task. For this purpose, by choosing $t > 1$ large enough to realize the relation $1 \leq \log t \leq t \leq t^{2}$ one has a decomposition of the desired integrand:
\[N_{1}^{(2)}(t) = \int_{\vert\xi\vert \geq \frac{\delta_{0}}{\sqrt{\log t}}}\frac{\sin^{2}(t\vert\xi\vert)}{\vert\xi\vert^{2}}\vert w_{1}(\xi)\vert^{2}d\xi\]
\[+ \int_{\frac{\delta_{0}}{\sqrt{\log t}} \geq \vert\xi\vert \geq \frac{\delta_{0}}{\sqrt{t}}}\frac{\sin^{2}(t\vert\xi\vert)}{\vert\xi\vert^{2}}\vert w_{1}(\xi)\vert^{2}d\xi + \int_{\frac{\delta_{0}}{\sqrt{t}} \geq \vert\xi\vert \geq \frac{\delta_{0}}{t}}\frac{\sin^{2}(t\vert\xi\vert)}{\vert\xi\vert^{2}}\vert w_{1}(\xi)\vert^{2}d\xi \]
\begin{equation}\label{i-22-10}
=: O_{1}(t) + O_{2}(t) + O_{3}(t).
\end{equation}
Let us estimate them in order.\\
First, one has
\[O_{1}(t) \leq \frac{\log t}{\delta_{0}^{2}}\int_{\vert\xi\vert \geq \frac{\delta_{0}}{\sqrt{\log t}}}\vert w_{1}(\xi)\vert^{2}d\xi\]
\begin{equation}\label{i-22-1}
\leq  \frac{\log t}{\delta_{0}^{2}}\Vert u_{1}\Vert^{2},\quad t \gg 1.
\end{equation}
For  $O_{2}(t)$ one can proceed as follows:
\[ O_{2}(t) \leq \Vert u_{1}\Vert_{L^{1}}^{2}\omega_{2}\int_{\frac{\delta_{0}}{\sqrt{t}}}^{\frac{\delta_{0}}{\sqrt{\log t}}}\frac{1}{r}dr\]
\begin{equation}\label{i-22-2}
\leq  \omega_{2}\Vert u_{1}\Vert_{L^{1}}^{2}2^{-1}\left(\log t - \log(\log t)\right),\quad t \gg 1.
\end{equation}
Lastly, we treat $O_{3}(t)$ to give the following estimate such that
\[O_{3}(t) \leq \Vert u_{1}\Vert_{L^{1}}^{2}\int_{\frac{\delta_{0}}{t} \leq \vert\xi\vert \leq \frac{\delta_{0}}{\sqrt{t}}}\frac{1}{r^{2}}d\xi\] 
\begin{equation}\label{i-22-3}
\leq  \omega_{2}\Vert u_{1}\Vert_{L^{1}}^{2}2^{-1}\log t,\quad t \gg 1.
\end{equation}
Thus, it follows from \eqref{i-22-10}, \eqref{i-22-1}, \eqref{i-22-2} and \eqref{i-22-3} that for large $t \gg 1$
\begin{equation}\label{i-22-4}
N_{1}^{(2)}(t) \leq C\left(\Vert u_{1}\Vert_{L^{1}}^{2} \log t + \Vert u_{1}\Vert^{2}\log t + \Vert u_{1}\Vert_{L^{1}}^{2}(\log t-\log(\log t))\right), \quad t \gg 1
\end{equation}
with some constant $C > 0$.
\noindent
Therefore, by combining \eqref{i-20}, \eqref{i-21} and \eqref{i-22-4} one can get the high-frequency estimate for $n = 2$:
\begin{equation}\label{i-high-2}
I_{high}^{(2)}(t) \leq C\left(\Vert u_{0}\Vert^{2} + \Vert u_{1}\Vert^{2} + \Vert u_{1}\Vert_{L^{1}}^{2}\right)\log t, \quad t \gg 1.
\end{equation}

Finally, it follows from \eqref{i-5}, \eqref{i-low}, \eqref{i-24} and \eqref{i-high-2} one can get the crucial upper bound estimates for each $n = 1$ and $n = 2$.
\begin{lem}\label{lem3}Let $n = 1$, and $[u_{0},u_{1}] \in (L^{1}({\bf R}^{n})\cap L^{2}({\bf R}^{n})) \times (L^{1}({\bf R}^{n})\cap L^{2}({\bf R}^{n}))$. Then, it holds that
\[\Vert w(t,\cdot)\Vert^{2} \leq C(\Vert u_{0}\Vert^{2} + \Vert u_{0}\Vert_{L^{1}}^{2} + \Vert u_{1}\Vert^{2} + \Vert u_{1}\Vert_{L^{1}}^{2})t, \quad t \gg 1.\]
\end{lem}
\begin{lem}\label{lem4}Let $n = 2$, and $[u_{0},u_{1}] \in (L^{1}({\bf R}^{n})\cap L^{2}({\bf R}^{n})) \times (L^{1}({\bf R}^{n})\cap L^{2}({\bf R}^{n}))$. Then, it holds that
\[\Vert w(t,\cdot)\Vert^{2} \leq C(\Vert u_{0}\Vert^{2} +\Vert u_{0}\Vert_{L^{1}}^{2} + \Vert u_{1}\Vert^{2} + \Vert u_{1}\Vert_{L^{1}}^{2})\log t, \quad t \gg 1.\]
\end{lem}
\par
\vspace{0.5cm}
Finally, the proofs of Theorems \ref{theorem1} and \ref{theorem2} are direct consequence of Lemmas \ref{lem1}, \ref{lem2}, \ref{lem3} and \ref{lem4}, and the Plancherel Theorem.\\


\section{Application to local energy decay problems}

In this section we apply Theorem \ref{theorem2} to the local energy decay problems. This problem is very important problem in the wave equation field, and one can observe many important previous papers (see \cite{AI, B, BK, Br, CI, D, I-1,I-2, IN, IS, LP, L, Mel, M, M-2, Mura, R, Ral, Secchi, S, ST, T, Tu, V, Vo, Wal, Z}), and a recent interesting application to inverse scattering problem of the local energy decay can be found in \cite{MS}. In particular, Morawetz \cite{M} studied such problems by constructing the so-called Morawetz identity to develop the multiplier method (see \cite{IK} for its text book). The author in \cite{M} considered the following exterior mixed problems: 
\begin{align}
& u_{tt} -\Delta u = 0,\ \ \ (t,x)\in (0,\infty)\times \Omega,\label{eqn-1}\\
& u(0,x)= u_0(x), \quad  u_{t}(0,x)= u_{1}(x),\ \ \ x\in \Omega,\label{initial-1}\\
& u(t,\sigma) = 0,\ \ \ \sigma \in \partial\Omega,
\end{align}
where $\Omega \subset {\bf R}^{n}$ is a smooth exterior domain with bounded boundary $\partial\Omega$ and star-shaped compliment $\Omega^{c}$. In \cite{M} the author derived the crucial estimate:
\[\int_{B_{R}\cap\Omega}\left(\vert u_{t}(t,x)\vert^{2} + \vert\nabla u(t,x)\vert^{2} \right)dx \leq C_{R}t^{-1},\quad t \gg 1,\] 
where $C_{R} > 0$ is a constant depending on $R > 0$ and initial data $[u_{0},u_{1}] \in H_{0}^{1}(\Omega)\times L^{2}(\Omega)$, and $B_{R} := \{x \in {\bf R}^{n}\,:\,\vert x\vert \leq R\}$ ($R > 0$). To study it, the author assumed the compactness of the support of the initial data (say {\bf CSI} in short). The most difficult point to derive the result is to show the $L^{2}$-boundedness of the solution, that is to say, we have to show that $\Vert u(t,\cdot)\Vert \leq \exists M$ for all $t \geq 0$. Such CSI was essentially used to derive the $L^{2}$-boundedness of the solution. So, a natural question arises: can one remove the CSI assumption? This question can be solved by Zachmanoglou \cite{Z}, Murave\u\i \cite{Mura}, and the Ikehata's series of previous papers \cite{AI, CI, I-1, I-2, IN, IS}. In particular, in \cite{IN} the authors developed a refined multiplier method based on \cite{IM} combined with the weighted energy method originally developed by Todorova-Yordanov \cite{TY} in order to accomplish the removal of CSI. In particular, the method introduced in \cite{IM} is useful to get such $L^{2}$-boundedness of the solution, and the essential idea in \cite{IM} is like this: if one has the so-called Hardy-type inequality, one can get the $L^{2}$-boundedness of the solution. In this connection, in \cite{M} the author never used the Hardy inequality. In stead of using the Hardy type one, in \cite{M} the author solved a corresponding Poisson equation to get $L^{2}$-boundedness of the solution using CSI. These methods can be also applied to the Cauchy problem (1.1)-(1.2) in ${\bf R}^{n}$, and in fact, without CSI condition Chara\~o-Ikehata \cite{CI} have derived the local energy decay results by using the method due to \cite{IN} to problem (1.1)-(1.2). However, it should be remarked that in \cite{CI} only the two dimensional case is not studied yet. 
This is because in the two dimensional case, one can not obtain the Hardy type inequality as was pointed out by \cite[(1.14) of Theorem 5]{O}. One of useful ideas to treat the two dimensional case is to impose the extra assumption $\int_{{\bf R}^{2}}u_{1}(x)dx = 0$ to get $L^{2}$-boundedness of the solution, and this is tried in \cite{Azer-ike} recently. The method in \cite{Azer-ike} is independent from establishing the Hardy-type inequality, and the Fourier transform is suitably used like in this paper. From this observation, in the whole space case one can apply the Fourier transform, and it suffices to analyze the $L^{2}$-behavior of the Fourier transformed solution $\hat{u}(t,\xi)$ through the Plancherel Theorem. At present, in the two dimensional case, under the non-vanishing condition of the moment $\int_{{\bf R}^{2}}u_{1}(x)dx$ it seems open to get the local energy decay because in general, the $L^{2}$-boundedness of the solution is not trivial.\\

Now, from Theorem \ref{theorem2} one can get the upper bound estimate with $log$-order under the loose condition on the initial data (not necessarily $\int_{{\bf R}^{2}}u_{1}(x)dx = 0$). Although this is not the boundedness, it may work nicely to get the local energy decay in the two dimensional whole space case. It should be mentioned once more that one does not rely on any Hardy type inequalities.\\     

Let us mention our story. First the Morawetz identity in ${\bf R}^{n}$-version can be stated as follows. Let $u \in C([0,\infty);H^{1}({\bf R}^{n})) \cap C^{1}([0,\infty);L^{2}({\bf R}^{n}))$ be the weak solution to problem (1.1)-(1.2).  Then, one has the following identity.
\[tE(t) = \frac{n-1}{2}\int_{{\bf R}^{n}}u_{1}(x)u_{0}(x)dx + \int_{{\bf R}^{n}}u_{1}(x)(x\cdot\nabla u_{0}(x))dx - \frac{n-1}{2}\int_{{\bf R}^{n}}u_{t}(t,x)u(t,x)dx\]
\begin{equation}\label{3-1}
-\int_{{\bf R}^{n}}u_{t}(t,x)(x\cdot\nabla u(t,x))dx,\quad t \geq 0.
\end{equation} 
$L^{2}$-boundedness directly affects the part $\int_{{\bf R}^{n}}u_{t}(t,x)u(t,x)dx$ because of the Schwarz inequality and energy identity (1.3):
\[\vert\int_{{\bf R}^{n}}u_{t}(t,x)u(t,x)dx\vert \leq \left(\int_{{\bf R}^{n}}\vert u(t,x)\vert^{2}dx\right)^{1/2}\left(\int_{{\bf R}^{n}}\vert u_{t}(t,x)\vert^{2}dx\right)^{1/2} = \Vert u(t,\cdot)\Vert\Vert u_{t}(t,\cdot)\Vert\]
\begin{equation}\label{3-2}
\leq \Vert u(t,\cdot)\Vert\sqrt{2E(0)},\quad t \geq 0.
\end{equation}
So, it suffices to control the quantity $\Vert u(t,\cdot)\Vert$ in order to get some good estimate of the term
$$F(t) := \int_{{\bf R}^{n}}u_{t}(t,x)u(t,x)dx.$$
\begin{rem}{\rm In the case when $n = 1$, the term $F(t)$ vanishes in \eqref{3-1}, so in the one dimensional case we only control the term
\[G(t) := \int_{{\bf R}^{n}}u_{t}(t,x)(x\cdot\nabla u(t,x))dx.\]
The analysis on $G(t)$ can be done in \cite{CI, IN, IS} by the (modified) weighted energy method due to \cite{TY} in order to absorb this term to a part of the total energy under the non-compact support assumption of the initial data. So, this topic on $G(t)$ is out of scope of this paper.}
\end{rem}

Now, one introduces the midway formula derived in \cite{CI, IN} to get the local energy decay in the case when the CSI condition is not assumed. \\
For each $R > 0$, we define the local energy in the region $B_{R} \subset {\bf R}^{2}$: 
\[E_{R}(t) := \int_{B_{R}}\left(\vert u_{t}(t,x)\vert^{2} + \vert\nabla u(t,x)\vert^{2}\right)dx.\]
Then, the result in \cite[(2.13) at p. 272]{IN} and/or \cite[(3.5) with $C(x) \equiv 1$ and $\eta = 0$]{CI} tell us the following fact. 
\begin{pro}\label{pro1}\,Let $n = 2$, $R > 0$, and assume that $[u_{0},u_{1}] \in H^{1}({\bf R}^{2})\times L^{2}({\bf R}^{2})$ further satisfies
\[\int_{{\bf R}^{2}}\vert x\vert\left( \vert u_{1}(x)\vert^{2} + \vert\nabla u_{0}(x)\vert^{2}\right)dx < +\infty.\] 
Then, the weak solution $u(t,x)$ to problem {\rm (1.1)-(1.2)} satisfies 
\[(t-R)E_{R}(t) \leq K_{0} + \frac{1}{2}\left\vert \int_{{\bf R}^{2}}u_{t}(t,x)u(t,x)dx\right\vert, \quad \forall t > R,\]
where
\[K_{0} := \int_{{\bf R}^{2}}u_{1}(x)(x\cdot\nabla u_{0}(x))dx + \frac{1}{2}\int_{{\bf R}^{2}}u_{1}(x)u_{0}(x)dx + E(0).\]
\end{pro}
Therefore, from this Proposition \ref{pro1} one can observe that if one can control the function $F(t)$ then one has a possibility to get the local energy decay. In \cite{CI} the authors have imposed rather stronger condition such that
\[\int_{{\bf R}^{2}}u_{1}(x)dx = 0\]
to get the "boundedness" of the function $F(t)$. This is their weak point in \cite{CI} in the two dimensional case. \\
Let us apply Theorem 2 to control \eqref{3-2}. Indeed, by applying Theorem 2 to \eqref{3-2} one can get  
\begin{equation}\label{3-3}
\left\vert F(t)\right\vert \leq C\sqrt{2E(0)}I_{0,2}\sqrt{\log t},\quad t \gg 1,
\end{equation}
where the constant $C > 0$ does not depend on any $R > 0$. By combining \eqref{3-3} with Proposition \ref{pro1} one has arrived at the following crucial result.
\begin{theo}\label{theorem5}\,Let $n = 2$, $R > 0$, and assume that $[u_{0},u_{1}] \in (H^{1}({\bf R}^{2})\cap L^{1}({\bf R}^{2}))\times(L^{2}({\bf R}^{2})\cap L^{1}({\bf R}^{2}))$ further satisfies
\[\int_{{\bf R}^{2}}\vert x\vert\left( \vert u_{1}(x)\vert^{2} + \vert\nabla u_{0}(x)\vert^{2}\right)dx < +\infty.\] 
Then, the weak solution $u(t,x)$ to problem {\rm (1.1)-(1.2)} satisfies
\[E_{R}(t) \leq \frac{K_{0}}{t-R} + \frac{C}{2}\sqrt{2E(0)}I_{0,2}\frac{\sqrt{\log t}}{t-R}, \quad \forall t > R.\]
\end{theo}
\begin{rem}{\rm It should be strongly mentioned that Theorem \ref{theorem5} can be derived without assuming CSI. For the one dimensional case, the effect of the function $F(t)$ is nothing, so one has no any good applications of Theorem 1.1 yet.}
\end{rem}
\begin{rem}
{\rm Shapiro \cite{S} has also obtained the local energy decay estimates for $n = 2$ with $(\log t)^{-2}$-order for more general wave equations with variable coefficients including (1.1), however, the work of \cite{S} has been done under stronger assumptions such as CSI and high regularity on the initial data. On the other hand, by \cite[Theorem 2.1]{Racke} applied to the case for $n = 2$ one can get the local energy decay with rate $t^{-1}$ under high regularity assumptions on the initial data (see also \cite[(22)]{Secchi}).}
\end{rem} 
\begin{rem}{\rm By combining our method with the Morawetz identity, in the forthcoming project we will study similar local energy decay estimates for the linear Klein-Gordon equation in the framework of non-compactly supported initial data.}
\end{rem}
\par
\vspace{0.5cm}
\noindent{\em Acknowledgement.}
\smallskip
The author would like to thank Masaru Ikehata (Hiroshima University) for his fruitful discussion at the preparation stage of this paper. 
The work of the author was supported in part by Grant-in-Aid for Scientific Research (C)20K03682  of JSPS.


\end{document}